\documentclass[a4paper,12pt]{article}
\usepackage{comment}
\usepackage{cite}
\usepackage{amsmath}
\usepackage{amssymb}
\usepackage{amsfonts}
\usepackage[T1]{fontenc}
\usepackage[utf8]{inputenc}
\usepackage{graphicx}
\usepackage{fancyhdr}
\usepackage{float}
\usepackage{xcolor}
\usepackage{authblk}
\usepackage{mathrsfs}
\usepackage{empheq}
\usepackage[hyphens]{url}
\usepackage{hyperref} 
\usepackage[]{breakurl}
\usepackage[]{amsthm}
\usepackage{tikz}

\usepackage{tikz}

\newcommand\catalannumber[3]{
  \fill[yellow]  (#1) rectangle +(#2,#2);
  \fill[fill=gray!15]
  (#1)
  \foreach \dir in {#3}{
    \ifnum\dir=0
    -- ++(1,0)
    \else
    -- ++(0,1)
    \fi
  } |- (#1);
  \draw[help lines] (#1) grid +(#2,#2);
  \draw[dashed] (#1) -- +(#2,#2);
  \coordinate (prev) at (#1);
  \foreach \dir in {#3}{
    \ifnum\dir=0
    \coordinate (dep) at (1,0);
    \else
    \coordinate (dep) at (0,1);
    \fi
    \draw[line width=2pt,-stealth] (prev) -- ++(dep) coordinate (prev);
  };
}

\newtheorem{theorem}{Theorem}

\usepackage{geometry}
\begin{document}

\title{An integral representation of Catalan numbers using the F\'eaux formula}

\author[$\dagger$]{Jean-Christophe {\sc Pain}$^{1,2,}$\footnote{jean-christophe.pain@cea.fr}\\
\small
$^1$CEA, DAM, DIF, F-91297 Arpajon, France\\
$^2$Universit\'e Paris-Saclay, CEA, Laboratoire Mati\`ere en Conditions Extr\^emes,\\ 
F-91680 Bruy\`eres-le-Ch\^atel, France
}

\maketitle

\begin{abstract}
We present an integral expression of the Catalan numbers, based on F\'eaux' integral representation of $\log\left[\Gamma(x)\right]$, $\Gamma$ being the usual Gamma function. The obtained formula may be the starting point of the derivation of new relations involving central binomial coefficients or Catalan numbers.
\end{abstract}

\section{Introduction}

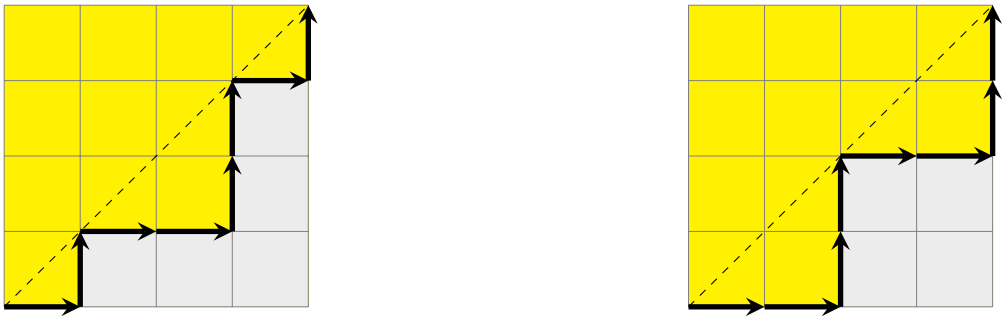
\begin{figure}
\centering
\begin{tikzpicture}
  \catalannumber{0,0}{4}{0,1,0,0,1,1,0,1}
  \catalannumber{9,0}{4}{0,0,1,1,0,0,1,1}
\end{tikzpicture}
\caption{Two monotonic lattice paths in the case $n=4$. They can be represented by listing the Catalan elements by column height: [0,1,1,3] (left) and [0,0,2,2] (right). There are 14 diagrams of the kinds of these two ones, which is precisely the fourth Catalan number $C_4$.}\label{fig1}
\end{figure}

The Catalan numbers are encountered in many fields of number theory and combinatorics \cite{Koshy2006,Koshy2008,Stanley2015,Larcombe1999}. They can be expressed in different ways, such as
\begin{equation*}
C_n=\frac{1}{n+1}\binom{2n}{n},
\end{equation*}
where $\binom{2n}{n}$ is the usual binomial coefficient, or also as
\begin{equation*}
C_n=~_2F_1\left[\begin{array}{cc}
1-n,-n\\
2
\end{array};1\right],
\end{equation*}
where $_2F_1$ represents the Gauss hypergeometric function. $C_n$ can be interpreted as the number of possibilities to split a convex polygon of $n+2$ sides into triangles by connecting vertices, or as the number of literal expressions containing $n$ pairs of correctly-matched parentheses. For instance, on the $n=4$ case, the different possibilities are ()()()(), ()()(()), ()(())(), ()(()()), ()((())), (())()(), (())(()), (()())(), ((()))(), (()()()), (()(())), ((())()), ((()())) and (((()))), their total number being $C_4=14$. As another example, $C_n$ is the number of Dyck words of length $2n$. A Dyck word is a string consisting of $n$ X and $n$ Y such that no segment of the string (starting from the full left) has more Y than X. The unique Dyck word of length 2 is XY ($C_1=1$), there are $C_2=2$ Dyck words of length 4: XXYY, XYXY, and $C_3=5$ Dyck words of length 6: XXXYYY, XYXXYY, XYXYXY, XXYYXY and XXYXYY. As a last example, it is worth mentioning that $C_n$ counts the monotonic lattice paths (i.e., starting from the lower left corner, consisting of right or up arrows and finishing in the upper right one) along the edges of a grid with $n \times n$ squares, which do not cross the diagonal \cite{Andre1887,Crepinsek2009} (two of them are displayed in Fig. \ref{fig1} in the case $n=4$). It turns out that numbering such paths is equivalent to counting the Dyck words (where X would stand for a right move and Y for an up move).

Many integral representations of Catalan numbers have been obtained \cite{Penson2001,Dana2012,Shi2015,Qi2017,Guo2023,Pain2024a,Pain2024b,Qi2024}. For instance, we recently obtained two representations \cite{Pain2024a} using the Touchard identity \cite{Touchard1928} on one hand, and a variant of the latter published by Callan \cite{Callan2013} on the other hand. We also found an expression using the Malmst\'en integral representation for the logarithm of the Gamma function \cite{Pain2024b}.

In section \ref{sec2} we present, using the F\'eaux integral (for the logarithm of the Gamma function), an integral representation for the Catalan numbers which, to our knowledge, was not published elsewhere. 

\section{Integral representation deduced from the F\'eaux formula}\label{sec2} 

\begin{theorem}

Let $n$ be a positive integer and $C_n$ the $n^{th}$ Catalan number. Then the integral representation
\begin{equation*}
    C_n=\frac{2^{2n}}{\sqrt{\pi}}\,\exp\left[\int_0^{\infty}\left\{\frac{\left[(1+t)^{3/2}-1\right]}{(1+t)^{n+2}\log(1+t)}-\frac{3}{2}e^{-t}\right\}\frac{\mathrm{d}t}{t}\right]
\end{equation*}
holds true.

\end{theorem}

\begin{proof}

Using the alternative representation of Catalan numbers:
\begin{equation*}
    C_n=\frac{4^n\Gamma(n+1/2)}{\sqrt{\pi}\,\Gamma(n+2)},
\end{equation*}
where
\begin{equation*}
\Gamma(x+1)=\int_0^{\infty}\,t^x\,e^{-t}\,\mathrm{dt}
\end{equation*}
is the usual Gamma function, one gets, for $n\geq 1$:
\begin{equation}\label{lncn}
    \log\left(C_n\right)=(2\log 2)\,n-\frac{\log\pi}{2} +\log\left[\Gamma\left(n+\frac{1}{2}\right)\right]-\log\left[\Gamma(n+2)\right].
\end{equation}
The F\'eaux formula reads \cite{Whittaker1990,Campbell1966,Gilbert1875}:
\begin{equation}\label{feaux}
    \log\left[\Gamma(x+1)\right]=\int_0^{\infty}\left[x\,e^{-t}+\frac{(1+t)^{-x-1}-(1+t)^{-1}}{\log(1+t)}\right]\frac{\mathrm{d}t}{t},
\end{equation}
yielding in particular, and for our purpose:
\begin{equation}\label{a1}
    \log\left[\Gamma\left(x+\frac{1}{2}\right)\right]=\int_0^{\infty}\left[\left(x-\frac{1}{2}\right)\,e^{-t}+\frac{(1+t)^{-x-1/2}-(1+t)^{-1}}{\log(1+t)}\right]\frac{\mathrm{d}t}{t},
\end{equation}
as well as
\begin{equation}\label{a2}
    \log\left[\Gamma(x+2)\right]=\int_0^{\infty}\left[(x+1)\,e^{-t}+\frac{(1+t)^{-x-2}-(1+t)^{-1}}{\log(1+t)}\right]\frac{\mathrm{d}t}{t}.
\end{equation}
Note that there is a typographical error in Ref. \cite{Campbell1966} (Ch. V, p. 189): $(1+t)$ should be replaced by $\log(1+t)$ in the denominator. Using Eqs. (\ref{a1}) and (\ref{a2}), one gets

\begin{equation*}
     \log\left[\Gamma\left(x+\frac{1}{2}\right)\right]-\log\left[\Gamma(x+2)\right]=\int_0^{\infty}\left\{\frac{\left[(1+t)^{3/2}-1\right]}{(1+t)^{n+2}\log(1+t)}-\frac{3}{2}e^{-t}\right\}\frac{\mathrm{d}t}{t},   
\end{equation*}
and Eq. (\ref{lncn}) becomes
\begin{equation*}
    \log\left(C_n\right)=(2\log 2)\,n-\frac{\log\pi}{2}+\int_0^{\infty}\left\{\frac{\left[(1+t)^{3/2}-1\right]}{(1+t)^{n+1}\log(1+t)}-\frac{3}{2}e^{-t}\right\}\frac{\mathrm{d}t}{t},   
\end{equation*}
and thus
\begin{equation*}
    C_n=\frac{2^{2n}}{\sqrt{\pi}}\,\exp\left[\int_0^{\infty}\left\{\frac{\left[(1+t)^{3/2}-1\right]}{(1+t)^{n+2}\log(1+t)}-\frac{3}{2}e^{-t}\right\}\frac{\mathrm{d}t}{t}\right],   
\end{equation*}
which completes the proof.

\end{proof}

\section{Tracks for deriving further new integral representations}

The Legendre duplication formula
\begin{equation}\label{dupli}
\Gamma (x)\,\Gamma\left(x+{\frac{1}{2}}\right)=2^{1-2x}\;{\sqrt {\pi }}\;\Gamma (2x),
\end{equation}
yields
\begin{equation*}
\log\left[\Gamma\left(x+\frac{1}{2}\right)\right]=(1-2x)\log 2+\frac{\log\pi}{2}+\log\left[\Gamma(2x)\right]-\log\left[\Gamma(x)\right].
\end{equation*}
Then, using Eq. (\ref{feaux}), one gets
\begin{equation*}
    \log\left[\Gamma(x)\right]=\int_0^{\infty}\left[(x-1)\,e^{-t}+\frac{(1+t)^{-x}-(1+t)^{-1}}{\log(1+t)}\right]\frac{\mathrm{d}t}{t},
\end{equation*}
as well as
\begin{equation*}
    \log\left[\Gamma(2x)\right]=\int_0^{\infty}\left[(2x-1)\,e^{-t}+\frac{(1+t)^{-2x}-(1+t)^{-1}}{\log(1+t)}\right]\frac{\mathrm{d}t}{t},
\end{equation*}
and thus
\begin{equation*}
\log\left[\Gamma\left(x+\frac{1}{2}\right)\right]=(1-2x)\,\log 2+\frac{\log\pi}{2}+\int_0^{\infty}\left[x\,e^{-t}+\frac{(1+t)^{-2x}-(1+t)^{-x}}{\log(1+t)}\right]\frac{\mathrm{d}t}{t},
\end{equation*}
which combined to Eq. (\ref{a2}) gives
\begin{align*}
\log\left[\Gamma\left(x+\frac{1}{2}\right)\right]-&\log\left[\Gamma(x+2)\right]=(1-2x)\,\log 2+\frac{\log\pi}{2}\nonumber\\
&+\int_0^{\infty}\left[\frac{(1+t)^{-2x}-(1+t)^{-x}-(1+t)^{-x-2}+(1+t)^{-1}}{\log(1+t)}-e^{-t}\right]\frac{\mathrm{d}t}{t},
\end{align*}
yielding an alternate integral representation of the Catalan numbers using Eq. (\ref{lncn}):
\begin{align*}
    \log(C_n)=&\log 2\nonumber\\
&+\int_0^{\infty}\left[\frac{(1+t)^{-2n}-(1+t)^{-n}-(1+t)^{-n-2}+(1+t)^{-1}}{\log(1+t)}-e^{-t}\right]\frac{\mathrm{d}t}{t},
\end{align*}
or
\begin{equation*}
    C_n=2\exp\left\{\int_0^{\infty}\left[\frac{(1+t)^{-2n}-(1+t)^{-n}-(1+t)^{-n-2}+(1+t)^{-1}}{\log(1+t)}-e^{-t}\right]\frac{\mathrm{d}t}{t}\right\}.
\end{equation*}

Alternate integral representation can also be derived combining, in Eqs. (\ref{dupli}) and (\ref{a2}), different integral representations for $\log\left[\Gamma\left(x\right)\right]$, $\log\left[\Gamma\left(x+1/2\right)\right]$, $\log\left[\Gamma\left(2x\right)\right]$ and  $\log\left[\Gamma\left(x+2\right)\right]$, such as the first Binet formula \cite{Sasvari1999}:
\begin{equation*}
\log\left[\Gamma(x)\right]=\left(x-{\frac {1}{2}}\right)\log x-x+{\frac {1}{2}}\log(2\pi)+\int_{0}^{\infty}\left({\frac{1}{2}}-\frac{1}{t}+\frac{1}{e^{t}-1}\right)\frac{e^{-tx}}{t}\,\mathrm{d}t,
\end{equation*}
the second Binet formula
\begin{equation*}
\log\left[\Gamma(x)\right]=\left(x-{\frac {1}{2}}\right)\log x-x+{\frac {1}{2}}\log(2\pi)+2\int_{0}^{\infty}{\frac{\arctan(t/x)}{e^{2\pi t}-1}}\,\mathrm{d}t,
\end{equation*}
the Malmst\'en integral \cite{Erdelyi1981,Pain2024b}:
\begin{equation*}
    \log\left[\Gamma(x)\right]=\int_0^{\infty}\left[(x-1)\,e^{-t}+\frac{(e^{-xt}-e^{-t})}{(1-e^{-t})}\right]\,\frac{\mathrm{d}t}{t}.
\end{equation*}
or the Kummer integral \cite{Whittaker1990,Kummer1847}:
\begin{equation*}
    \log\left[\Gamma(x)\right]=\frac{\log\pi}{2}-\frac{1}{2}\log\left[\sin(\pi x)\right]+\frac{1}{2}\int_0^{\infty}\left[\frac{\sinh\left[\left(\frac{1}{2}-x\right)t\right]}{\sinh\left(\frac{t}{2}\right)}-(1-2x)\,e^{-t}\right]\frac{\mathrm{d}t}{t}.
\end{equation*}
In addition, the Raabe formula:
\begin{equation*}
\int_{a}^{a+1}\log\left[\Gamma (x)\right]\,dx=\frac {1}{2}\log 2\pi +a\log a-a,
\end{equation*}
with $a>0$, may also lead to interesting relations (not necessarily integral representations, but rather sum rules for instance).

\section{Conclusion}\label{sec3}

In this work, an integral representation of the Catalan numbers was obtained, thanks to an integral of $\log\left[\Gamma(x+1)\right]$ due to F\'eaux, which is not as well-known as the Binet formulas for instance. In the future, we plan to apply similar techniques to the derivation of integral representations of the Ballot numbers
\begin{equation*}
B(n,k)=\frac{n-k}{n+k}\binom{n+k}{n} 
\end{equation*}
and of the Fuss-Catalan numbers:
\begin{equation*}
A_{m}(p,r)\equiv {\frac {r}{mp+r}}{\binom{mp+r}{m}}
\end{equation*}
as well as their multivariate extensions \cite{Aval2008}:
\begin{equation*}
    B_3(n,k,\ell)=\binom{n+k}{k}\binom{n+\ell-1}{\ell}\frac{n-k-\ell}{n+k}.
\end{equation*}

\end{document}